# Theoretical properties of Cook's PFC dimension reduction algorithm for linear regression

Oliver Johnson

*Statistics Group, Department of Mathematics, University of Bristol, University Walk, Bristol, BS8 1TW, UK*

**Abstract:** We analyse the properties of the Principal Fitted Components (PFC) algorithm proposed by Cook. We derive theoretical properties of the resulting estimators, including sufficient conditions under which they are $\sqrt{n}$-consistent, and explain some of the simulation results given in Cook's paper. We use techniques from random matrix theory and perturbation theory. We argue that, under Cook's model at least, the PFC algorithm should outperform the Principal Components algorithm.

## Contents

## 1. Introduction and notation

Cook [1] gives an extensive account of the history and controversy, dating back at least to the work of Fisher [4], surrounding the question of whether the selection of regression variables should be based only on values of the predictor, or whether the response should also be taken into account. For example, reduction by Principal Components (PC), introduced by Hotelling [7] and reviewed by for example Seber [12], only uses the sample covariance matrix of the predictors to select the variables, ignoring the values of the response.

Cook [1] argues that, while PC can play a useful role in regression problems, the analysis can also take account of the random response **Y**, and introduces the Principal Fitted Components (PFC) algorithm to do this. This algorithm performs a principal components analysis on the fitted sample covariance matrix obtained by projecting the predictor onto $\mathbf{f}_y$, a function of the response. Efficient selection of variables can allow regression algorithms to operate in a space of reduced dimension, making the resulting estimates potentially more accurate.

In this paper, we analyse the properties of the PFC algorithm, in the context of a model given by Equation (5) of Cook [1], following the notation of [1] throughout. Cook converts the more familiar forward linear model $\mathbf{Y} = \mathbf{X}\boldsymbol{\beta} + \boldsymbol{\varepsilon}$ into an equivalent inverse regression form:

**Example 1.1.** *Consider a sample of $n$ independent observations $\mathbf{X}_y$ in $\mathbb{R}^p$ indexed by the responses $y$ and generated as*

$$\mathbf{X}_y = \boldsymbol{\mu} + \boldsymbol{\Gamma}\boldsymbol{\beta}\mathbf{f}_y + \sigma\boldsymbol{\epsilon}. \tag{1}$$

*Here, $\boldsymbol{\mu}$ and each $\mathbf{X}_y$ are $p \times 1$ matrices (column vectors), $\boldsymbol{\Gamma}$ is a full rank $p \times d$ matrix (with $d < p$), $\boldsymbol{\beta}$ is a $d \times r$ matrix (with $d \leq r$) and $\mathbf{f}_y$ is an $r \times 1$ matrix. The parameter $\sigma$ gives the scale factor of $\boldsymbol{\epsilon}$, a $p \times 1$ matrix of errors. The entries of $\boldsymbol{\epsilon}$ will often (but not always) be assumed to be independent standard normals. We will estimate the span of the columns of the rank $d$ matrix $\boldsymbol{\Gamma}$, where by assumption, $\boldsymbol{\Gamma}^T\boldsymbol{\Gamma} = \mathbf{I}_d$ (we can reparameterize the model to achieve this). Having estimated $\boldsymbol{\Gamma}$, we can work in a space of smaller dimension, using whatever techniques are appropriate.*

The $\mathbf{f}_y$ is a vector-valued function of the random response $\mathbf{Y}$, and for simplicity we assume throughout that $\sum_y \mathbf{f}_y = \mathbf{0}$. Further we assume that $\boldsymbol{\epsilon}$ is independent of $\mathbf{Y}$, and hence of $\mathbf{f}_y$. (Note that we use different symbols, $\boldsymbol{\epsilon}$ and $\boldsymbol{\varepsilon}$, for the error terms in Equation (1) and the forward regression $\mathbf{Y} = \mathbf{X}\boldsymbol{\beta} + \boldsymbol{\varepsilon}$. We do not claim that $\boldsymbol{\varepsilon}$ is independent of $\mathbf{Y}$). The $\mathbf{f}_y$ can be constructed in a variety of ways, depending on the exact form of the data. For example, Cook mentions that if the conditional mean $\mathbb{E}(X|Y=y)$ can be modelled by a polynomial of degree $r$, then it is appropriate to take $\mathbf{f}_y = (y - m_1, y^2 - m_2, \ldots, y^r - m_r)^T$, where $m_u$ is the sample mean of the $u$th power of $y$. Alternatively, in the spirit of the Sliced Inverse Regression algorithm of Li [10], $\mathbf{f}_y$ can be constructed by slicing the range of $\mathbf{Y}$ into $(r+1)$ disjoint bins.

For ease of calculation we convert the model from Example 1.1 into matrix form. We write $\mathbb{X}$ for the $n \times p$ centred matrix of predictors, with rows $(\mathbf{X}_y - \overline{\mathbf{X}})^T$ (where $\overline{\mathbf{X}}$ is the sample mean of $\mathbf{X}$), write $\mathbf{F}$ for the $n \times r$ matrix with rows $\mathbf{f}_y^T$, and $\mathbf{E}$ for the $n \times p$ error matrix with rows $\sigma\boldsymbol{\epsilon}^T$. Given a full rank matrix $\mathbf{G}$ we write $\mathbf{P_G} = \mathbf{G}(\mathbf{G}^T\mathbf{G})^{-1}\mathbf{G}^T$ for the matrix which projects orthogonally onto the span of the columns of $\mathbf{G}$.

**Definition 1.2.** *Cook [1] defines the fitted matrix of predictors $\widehat{\mathbb{X}} = \mathbf{P_F}\mathbb{X}$ and proposes the PFC algorithm; an estimate $\widehat{\boldsymbol{\Gamma}}_{\mathrm{PFC}}$ of $\boldsymbol{\Gamma}$ is given by the set of $d$ eigenvectors of the fitted sample covariance matrix $\widehat{\mathbb{X}}^T\widehat{\mathbb{X}}$ which correspond to the largest eigenvalues.*

Cook contrasts this with the PC algorithm of Hotelling [7], which performs the corresponding calculation for the sample covariance matrix $\mathbb{X}^T\mathbb{X}$. That is, an estimate $\widehat{\boldsymbol{\Gamma}}_{\text{PC}}$ of $\boldsymbol{\Gamma}$ is given by the set of $d$ eigenvectors of $\mathbb{X}^T\mathbb{X}$ corresponding to the largest eigenvalues. Note that we refer to $\widehat{\boldsymbol{\Gamma}}_{\text{PFC}}$ and $\widehat{\boldsymbol{\Gamma}}_{\text{PC}}$ as estimates of $\boldsymbol{\Gamma}$ for the sake of brevity; in fact, the span of the columns of $\widehat{\boldsymbol{\Gamma}}_{\text{PFC}}$ and $\widehat{\boldsymbol{\Gamma}}_{\text{PC}}$ form estimates of the span of the columns of $\boldsymbol{\Gamma}$.

Theorem 2.4 of this paper gives a distributional result for the accuracy of PFC estimation. This allows us to explain the simulations relating to PFC estimation in Section 5 of Cook's paper [1], and to bound confidence intervals for the accuracy of the estimation, see Theorem 3.3 below. In Section 3 we give sufficient conditions for $\widehat{\boldsymbol{\Gamma}}_{\text{PFC}}$ to be a $\sqrt{n}$-consistent estimator of $\boldsymbol{\Gamma}$, in the case $d = r$. In Section 4, we use results from perturbation and random matrix theory to consider the more general case $d \leq r$ and consider the order of the errors that arise. In Theorem 4.3, we give sufficient conditions for $\widehat{\boldsymbol{\Gamma}}_{\text{PFC}}$ to be a $\sqrt{n}$-consistent estimator of $\boldsymbol{\Gamma}$ in this more general case.

Model (1) represents a problem in dimension reduction. The PC model is given as Equation (2) of [1] as

$$\mathbf{X}_y = \boldsymbol{\mu} + \boldsymbol{\Gamma}\boldsymbol{\nu}_y + \sigma\boldsymbol{\epsilon}. \qquad (2)$$

where $d \times 1$ vector $\boldsymbol{\nu}_y$ satisfies $\sum_y \boldsymbol{\nu}_y = 0$. Hence the model (1) is a special case of (2), where the $\boldsymbol{\beta}\mathbf{f}_y$ replaces $\boldsymbol{\nu}_y$. As Cook remarks, in the PFC model (1) we aim to estimate $\boldsymbol{\beta}$, which contains $dr$ parameters, whereas the PC model (2) contains the $(n-1)d$ parameters of $\boldsymbol{\nu}_y$. Hence, the PFC model has the attractive feature that the number of parameters to be estimated does not grow with $n$.

In Section 5 we argue that the PFC algorithm should perform strictly better than the PC algorithm, if the errors $\boldsymbol{\epsilon}$ are normally distributed. Specifically, Lemma 5.1 shows that in this case the sample covariance matrix $\mathbb{X}^T\mathbb{X}$ is the fitted covariance matrix $\widehat{\mathbb{X}}^T\widehat{\mathbb{X}}$ perturbed by random noise independent of $\widehat{\mathbb{X}}$ and $\boldsymbol{\Gamma}$. Hence inference about $\boldsymbol{\Gamma}$ using $\mathbb{X}$ (PC algorithm) will be necessarily less accurate than inference using $\widehat{\mathbb{X}}$ (PFC algorithm). Proposition 5.4 gives bounds in certain parameter regimes that imply that the PC estimate $\widehat{\boldsymbol{\Gamma}}_{\text{PC}}$ is also $O_{\mathbb{P}}(1/\sqrt{n})$, approximately explaining simulations relating to PC estimation in Section 5 of [1]. All the proofs of this paper are presented in Appendix A.

It is of course important to consider the validity of the model given by Equation (1), in order to understand the significance of these results. The key is the reversal of the conditional distribution for $\mathbf{Y}|\mathbf{X}$ implied by the linear model $\mathbf{Y} = \mathbf{X}\boldsymbol{\beta} + \boldsymbol{\varepsilon}$ to give a conditional distribution for $\mathbf{X}|(\mathbf{Y} = \mathbf{y})$, as in Equation (2). Such a reversal is not new, arising for example in the work of Oman [11]. If the $\mathbf{X}$ and $\boldsymbol{\varepsilon}$ are each multivariate normal, then $(\mathbf{X}, \mathbf{Y})$ will be jointly multivariate normal, and we can move between $\mathbf{Y} = \mathbf{X}\boldsymbol{\beta} + \boldsymbol{\varepsilon}$ and Equation (2) by parameterizing appropriately. The question of the validity of the full PFC model Equation (1) amounts to the issue of how accurately the $\boldsymbol{\beta}\mathbf{f}_y$ models $\boldsymbol{\nu}_y$. Presumably for real data there will be a trade-off between improved accuracy arising from PFC estimation and errors introduced by $\boldsymbol{\beta}\mathbf{f}_y$ not adequately explaining $\boldsymbol{\nu}_y$. Implicitly, model (2) requires that the $\mathbf{X}$ should be random and take continuous values. Hence $\mathbf{X}$ should come from measurements rather than designed experiments, and the case of factor variables is excluded.

In future work, we hope to explain the performance (in the sense of Mean Squared Error) of prediction errors shown in Figure 1(d) of Cook's paper [1], and to investigate the theoretical properties of the PFC$_{\text{PC}}$ algorithm described in Section 6 of [1].

## 2. Theoretical PFC performance when $d = r$

In this section, we analyse the performance of the PFC algorithm in the case $d = r$. In Lemma 2.2 we give explicit expressions for the matrices involved, and define a matrix $\mathbf{V}$ which gives the span of the PFC directions. In Lemma 2.3, we deduce distributional results for $\mathbf{V}$, which we use to prove Theorem 2.4. Finally, we show how Theorem 2.4 explains some of the simulation results given in Section 5 of [1].

To measure the accuracy of our estimates, we consider the distribution of a normalised version of the quantity $m(\widehat{\mathbf{\Gamma}}, \mathbf{\Gamma}) = \|(\mathbf{I}_p - \mathbf{P}_{\mathbf{\Gamma}})\widehat{\mathbf{\Gamma}}\|$ used, for example, by Xia, Tong, Li and Zhu [15]. We have a choice of which matrix norm $\|\cdot\|$ to use; we shall use the Frobenius norm $\|\cdot\|_F$ defined by

$$\|\mathbf{A}\|_F^2 = \text{tr}\,(\mathbf{A}\mathbf{A}^T) = \text{tr}\,(\mathbf{A}^T\mathbf{A}) = \sum_k \left|\mathbf{A}^{(k)}\right|^2,$$

where $\mathbf{A}^{(k)}$ is the $k$th column of $\mathbf{A}$ and $|\cdot|$ represents the vector norm. This choice of matrix norm has the attractive feature that $\|\mathbf{A}\|_F = \|\mathbf{AP}\|_F$ for orthogonal $\mathbf{P}$, so that we can make orthogonal changes of basis without affecting $m(\widehat{\mathbf{\Gamma}}, \mathbf{\Gamma})$.

**Definition 2.1.** *For true $\mathbf{\Gamma}$ and estimated value $\widehat{\mathbf{\Gamma}}$, we define*

$$C(\widehat{\mathbf{\Gamma}}, \mathbf{\Gamma}) = \frac{\|\mathbf{P}_{\mathbf{\Gamma}}\widehat{\mathbf{\Gamma}}\|_F^2}{\|(\mathbf{I}_p - \mathbf{P}_{\mathbf{\Gamma}})\widehat{\mathbf{\Gamma}}\|_F^2}.$$

In the case $d > 1$, the quantity $C(\widehat{\mathbf{\Gamma}}, \mathbf{\Gamma})$ will change on rescaling columns of the matrix $\widehat{\mathbf{\Gamma}}$ by different amounts (that is, on multiplying $\widehat{\mathbf{\Gamma}}$ by diagonal $\mathbf{D}$ not proportional to $\mathbf{I}$). We do not have a priori estimates for the true scalings of the eigenvectors, so use the scalings that arise from orthogonal transformations of the columns of a matrix $\mathbf{V}$ arising from matrix factorization (see Lemma 2.2). Lemma 2.3 shows that this choice has the attractive feature that the entries of $\mathbf{V}$ have equal variance. There may exist better scalings in the sense of reduced average angle, however this choice already has good properties, as the theorems of this paper show.

**Lemma 2.2.** *Under model (1) from Example 1.1, the span of the $r$ largest eigenvectors of $\widehat{\mathbb{X}}^T\widehat{\mathbb{X}}$ is identical to the span of columns of the $p \times r$ matrix $\mathbf{V}$ defined by*

$$\mathbf{V} = \mathbf{\Gamma}\boldsymbol{\beta}(\mathbf{F}^T\mathbf{F})^{1/2} + \mathbf{E}^T\mathbf{F}(\mathbf{F}^T\mathbf{F})^{-1/2}. \tag{3}$$

*In particular, in the case where $d = r$ these columns define the PFC space. Note that $\mathbf{V}$ can be found from known information and without calculating eigenvalues, by $\mathbf{V}^T = (\mathbf{F}^T\mathbf{F})^{-1/2}\mathbf{F}^T\widehat{\mathbb{X}} = (\mathbf{F}^T\mathbf{F})^{-1/2}\mathbf{F}^T\mathbb{X}.$*

The proofs of all the results of this paper are given in Appendix A. Next we give distributional results for $\mathbf{V}$, by conditioning on the values of $\mathbf{F}$.

**Lemma 2.3.** *Assume that the errors $\boldsymbol{\epsilon}$ have mean zero and are uncorrelated with variance 1. Conditioned on the values of $\mathbf{F}$, the entries of $\mathbf{V}$ have mean $\boldsymbol{\Gamma}\boldsymbol{\beta}(\mathbf{F}^T\mathbf{F})^{1/2}$ and variance $\sigma^2$, and are uncorrelated.*

Lemma 2.3 implies that if the errors $\boldsymbol{\epsilon}$ are independent standard normals, then the columns of $\mathbf{V}$ are independent multivariate normals, with means given by the columns of $\boldsymbol{\Gamma}\boldsymbol{\beta}(\mathbf{F}^T\mathbf{F})^{1/2}$ and with covariance matrix $\sigma^2\mathbf{I}_p$. Hence the $p \times p$ matrix $\widehat{\mathbb{X}}^T\widehat{\mathbb{X}} = \mathbf{V}\mathbf{V}^T$ has a non-central Wishart distribution with $r$ degrees of freedom, scale parameter $\sigma^2\mathbf{I}_p$ and non-centrality parameter $\boldsymbol{\Gamma}\boldsymbol{\beta}(\mathbf{F}^T\mathbf{F})\boldsymbol{\beta}^T\boldsymbol{\Gamma}^T$. This allows us to deduce a distributional result for $C$.

**Theorem 2.4.** *Under the model given by Equation (1), assume that the errors $\boldsymbol{\epsilon}$ are independent standard normals, and consider the PFC estimator $\widehat{\boldsymbol{\Gamma}}_{\mathrm{PFC}}$. In the case $d = r$, conditional on $(\mathbf{F}^T\mathbf{F})$, the term*

$$\frac{(p-d)C(\widehat{\boldsymbol{\Gamma}}_{\mathrm{PFC}},\boldsymbol{\Gamma})}{d} = \frac{(p-d)\|\mathbf{P}_{\boldsymbol{\Gamma}}\widehat{\boldsymbol{\Gamma}}_{\mathrm{PFC}}\|_F^2}{d\|(\mathbf{I}_p-\mathbf{P}_{\boldsymbol{\Gamma}})\widehat{\boldsymbol{\Gamma}}_{\mathrm{PFC}}\|_F^2} \sim F_{rd,r(p-d)}(\Lambda),$$

*where $F_{rd,r(p-d)}(\Lambda)$ denotes a non-central F distribution with $(rd, r(p-d))$ degrees of freedom and non-centrality parameter given by the scaled trace $\Lambda = \mathrm{tr}\,(\boldsymbol{\beta}(\mathbf{F}^T\mathbf{F})\boldsymbol{\beta}^T)/\sigma^2$.*

The quantity $C(\widehat{\boldsymbol{\Gamma}},\boldsymbol{\Gamma})$ measures the proportion of the magnitude of the estimate $\widehat{\boldsymbol{\Gamma}}$ which lies in the span of the columns of $\boldsymbol{\Gamma}$, and hence measures how good an estimate of the span of $\boldsymbol{\Gamma}$ is provided by $\widehat{\boldsymbol{\Gamma}}$. In the case $r = d = 1$, this is compatible with Cook's plots of the angle $\Theta(\widehat{\boldsymbol{\Gamma}},\boldsymbol{\Gamma})$ between true $\boldsymbol{\Gamma}$ and estimated $\widehat{\boldsymbol{\Gamma}}$, in the sense that for any $\widehat{\boldsymbol{\Gamma}}$ and $\boldsymbol{\Gamma}$ the $C(\widehat{\boldsymbol{\Gamma}},\boldsymbol{\Gamma}) = \cot^2\Theta(\widehat{\boldsymbol{\Gamma}},\boldsymbol{\Gamma})$.

Section 6 of Li, Zha and Chiaromonte [9] introduces a different measure of similarity of subspaces as $\|\mathbf{P}_{\boldsymbol{\Gamma}} - \mathbf{P}_{\widehat{\boldsymbol{\Gamma}}}\|$, where $\|\cdot\|$ represents the operator norm. In the case $d = r = 1$, this again is compatible with the angle $\Theta(\widehat{\boldsymbol{\Gamma}},\boldsymbol{\Gamma})$, in the sense that in this case $\|\mathbf{P}_{\boldsymbol{\Gamma}} - \mathbf{P}_{\widehat{\boldsymbol{\Gamma}}}\| = |\sin\Theta(\widehat{\boldsymbol{\Gamma}},\boldsymbol{\Gamma})|$. The use of the operator norm means that the measure of [9] represents 'worst case' performance, whereas using the Frobenius norm gives 'average' performance. Our techniques do not at present give distributional results for $\|\mathbf{P}_{\boldsymbol{\Gamma}} - \mathbf{P}_{\widehat{\boldsymbol{\Gamma}}}\|$, however, Frobenius norms are typically easier to calculate than operator norms.

Theorem 2.4 shows that the distribution of $C(\widehat{\boldsymbol{\Gamma}}_{\mathrm{PFC}},\boldsymbol{\Gamma})$ does not depend on the value of $\boldsymbol{\Gamma}$ itself, helping to explain Cook's remark [1, P.10] that "the value of principal component estimators does not rest solely with the presence of collinearity". Further, using Theorem 2.4, we can better understand the simulation graphs given in Section 5 of Cook [1].

**Example 2.5.** *Figure 1 of [1] considers the model given by Equation (1) in the case where $p = 10$, $d = r = 1$, $\boldsymbol{\beta} = 1$. Further, the errors $\boldsymbol{\epsilon}$ are normally distributed and $\mathbf{Y}$ is normal with variance $\sigma_Y^2$, so that $\mathbf{F}^T\mathbf{F} = \sum_{i=1}^n (y_i -$* $\overline{y})^2 \sim \sigma_Y^2 \chi_{n-1}^2$. *We consider the angle* $\Theta(\widehat{\boldsymbol{\Gamma}}_{\mathrm{PFC}}, \boldsymbol{\Gamma})$, *where* $\cot^2 \Theta(\widehat{\boldsymbol{\Gamma}}_{\mathrm{PFC}}, \boldsymbol{\Gamma}) = C(\widehat{\boldsymbol{\Gamma}}_{\mathrm{PFC}}, \boldsymbol{\Gamma})$.

*In Figure 1 we simulate directly from the distribution given by Theorem 2.4 and vary parameters* $n$, $\sigma$ *and* $\sigma_Y$. *Based on a sample of size 50,000 for each set of parameter values, we plot the sample mean as* ∘, *along with upper and lower 5% sample quantiles (plotted as* △ *and* +*). The sample means shown here in Figure 1(a)-1(c) fit closely with those in Figure 1(a)–1(c) of [1].*

## 3. $\sqrt{n}$-consistency of PFC estimates

The quantiles plotted in Figure 1 give some idea of the spread of likely values of the angle $\Theta(\widehat{\boldsymbol{\Gamma}}_{\mathrm{PFC}}, \boldsymbol{\Gamma})$. Theorem 2.4 shows that the squared cotangent $C(\widehat{\boldsymbol{\Gamma}}_{\mathrm{PFC}}, \boldsymbol{\Gamma})$ is a scalar multiple of a non-central $F$ distribution with random non-centrality parameter, a complicated hierarchical form of mixture distribution, meaning that it is not trivial to give confidence intervals for the angle $\Theta(\widehat{\boldsymbol{\Gamma}}_{\mathrm{PFC}}, \boldsymbol{\Gamma})$ in closed form.

In the case $d = r$ we give probabilistic bounds in Theorem 3.3 demonstrating that $\Theta(\widehat{\boldsymbol{\Gamma}}_{\mathrm{PFC}}, \boldsymbol{\Gamma})$ decays like $1/\sqrt{n}$ (so $\widehat{\boldsymbol{\Gamma}}_{\mathrm{PFC}}$ is a $\sqrt{n}$-consistent estimator of $\boldsymbol{\Gamma}$), for a more general class of error models than simply assuming normality. For simplicity of exposition, we restrict to the case where the distributions of $\boldsymbol{\epsilon}$ are symmetric, though this is not necessary for $\sqrt{n}$-consistency to hold – see Appendix A for details. First we give two technical results, Lemmas 3.1 and 3.2, concerning the entries of the matrix $\mathbf{V}$ introduced in Lemma 2.2 and matrix $\boldsymbol{\beta}(\mathbf{F}^T\mathbf{F})\boldsymbol{\beta}^T$ respectively.

**Lemma 3.1.** *If the errors* $\boldsymbol{\epsilon}$ *are independent and symmetric with variance 1 and finite 4th moment* $m_4$ *then writing* $N = \|(\mathbf{I} - \mathbf{P}_{\boldsymbol{\Gamma}})\mathbf{V}\|_F^2$ *we know that*

$$\mathbb{E} N = \mathrm{tr}\,(\boldsymbol{\beta}(\mathbf{F}^T\mathbf{F})\boldsymbol{\beta}^T) + rd\sigma^2 \quad and \quad \mathrm{Var}\,(N) \leq dT + 4\sigma^2 \mathrm{tr}\,(\boldsymbol{\beta}(\mathbf{F}^T\mathbf{F})\boldsymbol{\beta}^T), \quad (4)$$

*where* $T = \sigma^4 r(m_4 - 1)$ *does not depend on* $n$. *Similarly, writing* $D = \|\mathbf{P}_{\boldsymbol{\Gamma}}\mathbf{V}\|_F^2$ *we know that under the same conditions*

$$\mathbb{E} D = r(p-d)\sigma^2 \quad and \quad \mathrm{Var}\,(D) \leq (p-d)T. \quad (5)$$

In the case where the errors $\boldsymbol{\epsilon}$ are standard normal, Lemma 3.1 simplifies, since the $\mathbf{C}_{ik}$ become independent and identically distributed normals. Thus $D/\sigma^2$ is central $\chi^2$ with $r(p-d)$ degrees of freedom. Similarly $N/\sigma^2$ is non-central $\chi^2$ with $rd$ degrees of freedom and non-centrality parameter $\mathrm{tr}\,(\boldsymbol{\beta}(\mathbf{F}^T\mathbf{F})\boldsymbol{\beta}^T)/\sigma^2$, and $T = 2r\sigma^4$, with equality holding in the bounds on $\mathrm{Var}\,(N)$ and $\mathrm{Var}\,(D)$ in Equations (4) and (5).

We will write $\lambda_1(\mathbf{X}) \geq \lambda_2(\mathbf{X}) \geq \ldots \geq \lambda_p(\mathbf{X})$ for the ordered sequence of eigenvalues of a real symmetric $p \times p$ matrix $\mathbf{X}$.

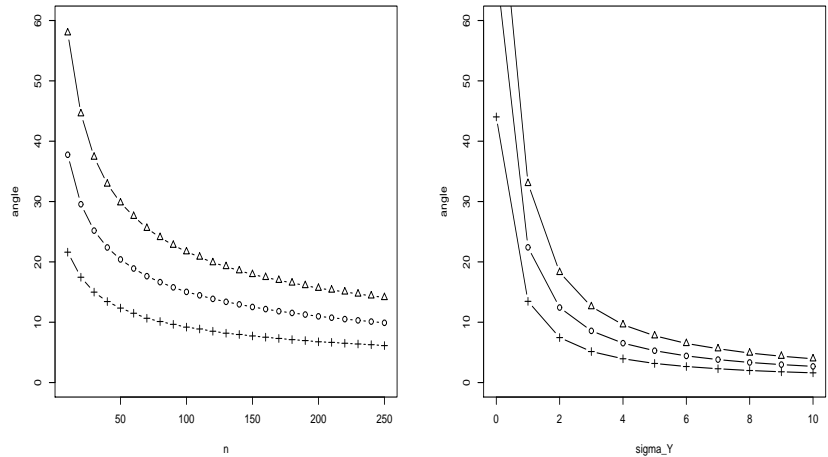

(a) Angle $\Theta(\widehat{\mathbf{\Gamma}}_{\mathrm{PFC}}, \mathbf{\Gamma})$ vs $n$

(b) Angle $\Theta(\widehat{\mathbf{\Gamma}}_{\mathrm{PFC}}, \mathbf{\Gamma})$ vs $\sigma_Y$

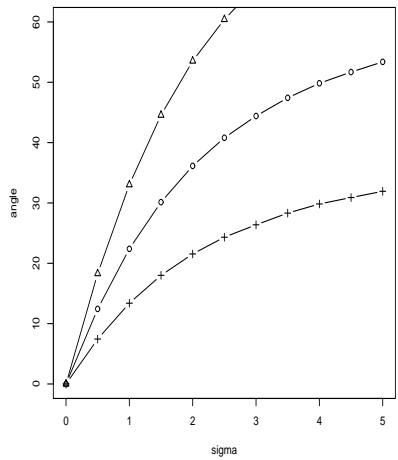

(c) Angle $\Theta(\widehat{\mathbf{\Gamma}}_{\mathrm{PFC}}, \mathbf{\Gamma})$ vs $\sigma$

Fig 1. *Simulation of angle $\Theta(\widehat{\mathbf{\Gamma}}_{\mathrm{PFC}}, \mathbf{\Gamma})$ between $\widehat{\mathbf{\Gamma}}_{\mathrm{PFC}}$ and $\mathbf{\Gamma}$, in setting of Example 2.5. We plot the sample mean as ○, and the upper and lower 5% sample quantiles as △ and +. (a) Angle vs n, with $\sigma_Y = \sigma = 1$; (b) Angle vs $\sigma_Y$, with $n = 40$, $\sigma = 1$; (c) Angle vs $\sigma$, with $n = 40$, $\sigma_Y = 1$.*

**Lemma 3.2.** *There exists a sequence* $(\phi_i : 1 \leq i \leq d)$ *such that given* $\epsilon > 0$ *and* $\delta > 0$, *there exists* $n^* = n^*(\delta, \epsilon)$ *such that the eigenvalues satisfy*

$$\mathbb{P}\left(\phi_i - \delta \leq \frac{\lambda_i(\boldsymbol{\beta}(\mathbf{F}^T\mathbf{F})\boldsymbol{\beta}^T)}{n} \leq \phi_i + \delta \text{ for all } i = 1, \ldots, d, \text{ and } n \geq n^*\right) \geq 1 - \epsilon.$$

**Theorem 3.3.** *In the case where* $d = r$ *and errors* $\boldsymbol{\epsilon}$ *are independent and symmetric with variance 1 and finite 4th moment, then we can construct confidence intervals such that*

$$\mathbb{P}\left(\Theta(\widehat{\boldsymbol{\Gamma}}_{\mathrm{PFC}}, \boldsymbol{\Gamma}) \geq \Theta^*_+(\alpha)\right) \leq \alpha \quad \text{and} \quad \mathbb{P}\left(\Theta(\widehat{\boldsymbol{\Gamma}}_{\mathrm{PFC}}, \boldsymbol{\Gamma}) \leq \Theta^*_-(\alpha)\right) \leq \alpha,$$

*where for any fixed* $\alpha$, *the* $\Theta^*_\pm(\alpha) = O(1/\sqrt{n})$.

Using Theorem 3.3, we can consider the case analysed in Example 2.5 and Figure 1, where $\boldsymbol{\beta} = 1$, $r = d = 1$, $p = 10$ and the errors $\boldsymbol{\epsilon}$ and $Y$ are normal. In this case, $(\mathbf{F}^T\mathbf{F}) \sim \sigma_Y^2 \chi_{n-1}^2$, so we can take $K_1 = \sigma_Y^2 - \epsilon$ and $K_2 = \sigma_Y^2 + \epsilon$. Equations (26) and (27) below show that confidence intervals for the angle $\Theta(\widehat{\boldsymbol{\Gamma}}_{\mathrm{PFC}}, \boldsymbol{\Gamma})$ decay asymptotically as $c/\sqrt{n}$, $c/\sigma_Y$ or $c\sigma$ respectively, as the other terms are kept constant, as Figure 1 may suggest.

## 4. Random matrices and perturbation

In this section, we analyse the general case $d \leq r$ of the model given by Equation (1) under the assumption that the errors $\boldsymbol{\epsilon}$ are Gaussian, using a perturbation argument. We first review some results from random matrix theory, which we use to deduce the order of $\Theta(\widehat{\boldsymbol{\Gamma}}_{\mathrm{PFC}}, \boldsymbol{\Gamma})$ in Theorem 4.3.

**Proposition 4.1.** *Write* $\mathbf{X}_{u,v}$ *for a* $u \times v$ *matrix with entries that are independent standard Gaussians, and consider the largest eigenvalue of the Wishart matrix* $\mathbf{X}_{u,v}\mathbf{X}_{u,v}^T$.

1. *For a sequence of matrices* $\mathbf{X}_{u_v,v}$, *in the regime where* $v \to \infty$ *and* $u_v/v \to \beta$,

$$\frac{1}{v}\lambda_1(\mathbf{X}_{u_v,v}\mathbf{X}_{u_v,v}^T) \longrightarrow (1 + \sqrt{\beta})^2 \quad \text{almost surely.} \tag{6}$$

2. *In the same regime*

$$\frac{\lambda_1(\mathbf{X}_{u_v,v}\mathbf{X}_{u_v,v}^T) - (\sqrt{v-1} + \sqrt{u_v})^2}{\sigma_{u_v v}} \xrightarrow{\mathcal{D}} F_1, \tag{7}$$

*where* $\sigma_{uv} = (\sqrt{v-1} + \sqrt{u})(1/\sqrt{v-1} + 1/\sqrt{u})^{1/3}$, *and where* $F_1$ *is the Tracy–Widom law of order 1.*

3. *For any* $u$ *and* $v$ *and for any* $t > 0$,

$$\mathbb{P}\left(\lambda_1(\mathbf{X}_{u,v}\mathbf{X}_{u,v}^T) \geq (\sqrt{u} + \sqrt{v} + t)^2\right) < \exp(-t^2/2). \tag{8}$$

Part 1 of Proposition 4.1 is due to Geman [6], who does not require that the entries of the matrix $\mathbf{X}_{u,v}$ be Gaussian (finiteness of moments of all orders will suffice). However, Equation (8) does require Gaussian entries, which is the reason for the assumption that $\boldsymbol{\epsilon}$ has Gaussian entries in the rest of this paper. Johnstone [8] proved Part 2, giving simulation results suggesting that this approximation is accurate for $u_v, v$ as small as 5. Part 3 is implied by Theorem II.13 of Davidson and Szarek [3].

We combine the random matrix results of Proposition 4.1 with a perturbation argument based on that given by Sibson [13] and used by Critchley [2] in a related context. As in Sibson [13], the perturbed eigenvectors are given in terms of generalized inverse matrices $\mathbf{M}^+$. If $\mathbf{Mx} = \lambda \mathbf{x}$, the linear map $\mathbf{M}^+$ is defined using the property that

$$\mathbf{M}^+ \mathbf{x} = \begin{cases} \lambda^{-1} \mathbf{x} & \text{if } \lambda \neq 0, \\ \mathbf{0} & \text{if } \lambda = 0. \end{cases}$$

To motivate the proof of $\sqrt{n}$-consistency in Theorem 4.3, we change basis to the orthogonal set $\{\mathbf{b}^{(i)}\}$ used in the proof of Theorem 2.4. Equations (17) and (19) below mean that in this new basis, $\widehat{\mathbb{X}}^T \widehat{\mathbb{X}} = (\mathbf{U} + \sigma \mathbf{S})(\mathbf{U} + \sigma \mathbf{S})^T$. Here $\mathbf{U}\mathbf{U}^T$ has a $d \times d$ block of the form $\boldsymbol{\beta}(\mathbf{F}^T\mathbf{F})\boldsymbol{\beta}^T$ with the remaining entries being zero, and $\mathbf{S}$ is a $p \times r$ matrix whose entries are independent standard Gaussians.

Hence if $\sigma = 0$ then $\widehat{\mathbb{X}}^T \widehat{\mathbb{X}}$ has $d$ positive eigenvalues with eigenvectors lying in the space spanned by the columns of $\boldsymbol{\Gamma}$, and the remaining eigenvalues are zero. Using perturbation theory, we bound how large $\sigma$ would have to be before one of the zero eigenvalues could become one of the $d$ largest ones.

**Definition 4.2.** *For $\sigma \neq 0$, we take $\mathbf{B} = \mathbf{U}\mathbf{U}^T$ and $\mathbf{L} = \sigma(\mathbf{S}\mathbf{U}^T + \mathbf{U}\mathbf{S}^T) + \sigma^2 \mathbf{S}\mathbf{S}^T$, so that $\widehat{\mathbb{X}}^T\widehat{\mathbb{X}} = \mathbf{B} + \mathbf{L}$. We define the first level crossing event*

$$\mathcal{L}_1 = \{\lambda_1(\mathbf{L}) - \lambda_p(\mathbf{L}) \geq \lambda_d(\boldsymbol{\beta}(\mathbf{F}^T\mathbf{F})\boldsymbol{\beta}^T)\}.$$

Since $\lambda_i(\mathbf{U}\mathbf{U}^T) = 0$ for $i \geq d+1$, results of Weyl (see Lemma A.1 below) imply that

$$\begin{aligned}\lambda_i(\widehat{\mathbb{X}}^T\widehat{\mathbb{X}}) &\geq \lambda_d(\boldsymbol{\beta}(\mathbf{F}^T\mathbf{F})\boldsymbol{\beta}^T) + \lambda_p(\mathbf{L}) & \text{for } i \leq d, \\ \lambda_i(\widehat{\mathbb{X}}^T\widehat{\mathbb{X}}) &\leq \lambda_1(\mathbf{L}) & \text{for } i \geq d+1.\end{aligned}$$

Hence, if $\mathcal{L}_1$ does not occur, the $d$ largest eigenvalues of $\widehat{\mathbb{X}}^T\widehat{\mathbb{X}}$ correspond to the perturbed values of the original eigenvalues of $\boldsymbol{\beta}(\mathbf{F}^T\mathbf{F})\boldsymbol{\beta}^T$. Proposition 4.1 gives probabilistic bounds on $\lambda_1(\mathbf{L}) - \lambda_p(\mathbf{L})$, allowing us to control $\mathbb{P}(\mathcal{L}_1)$.

**Theorem 4.3.** *In the case $d \leq r$, if the errors $\boldsymbol{\epsilon}$ are independent standard normals and the limiting matrix $\Phi = \lim_{n \to \infty}(\boldsymbol{\beta}(\mathbf{F}^T\mathbf{F})\boldsymbol{\beta}^T)/n$ has distinct eigenvalues, then there exist confidence intervals*

$$\mathbb{P}\left(\Theta(\widehat{\boldsymbol{\Gamma}}_{\text{PFC}}, \boldsymbol{\Gamma}) \geq \Theta^*(\alpha)\right) \leq \alpha,$$

*where for any fixed $\alpha$, the $\Theta^*(\alpha) = O(1/\sqrt{n})$.*

## 5. Theoretical performance of the PC algorithm

Finally we discuss the PC algorithm, arguing in Lemma 5.1 that, in the case of normal errors $\boldsymbol{\epsilon}$, it will give inferior performance to the PFC algorithm. Proposition 5.4 gives bounds on the performance of the PC algorithm, explaining the simulation results of Cook [1].

**Lemma 5.1.** *Under the model given by Equation (1), if the errors $\boldsymbol{\epsilon}$ are normally distributed, then we can write $\mathbb{X}^T \mathbb{X} = \widehat{\mathbb{X}}^T \widehat{\mathbb{X}} + N$, where $N$ is independent of $\widehat{\mathbb{X}}$ and $\boldsymbol{\Gamma}$.*

This result allows us to argue that PFC estimation should out-perform PC estimation, in several senses. Firstly, inference about random variable $X$ through random variable $Y$ is better (in the sense of Minimum Mean Squared Error) than inference about $X$ through $Y + N$, if $N$ is independent of $X$. This follows since the best estimates are $f(Y) = \mathbb{E}(X|Y)$ and $g(Y + N) = \mathbb{E}(X|Y + N)$ respectively. The fact that the MMSE is lower in the first case is equivalent to the fact that $\mathbb{E}g(Y + N)^2 \leq \mathbb{E}f(Y)^2$, which follows by the conditional Jensen inequality. Similarly, the conditional entropy $H(X|Y) \leq H(X|Y + N)$, showing that there is less uncertainty about $X$ on learning $Y$ than on learning $Y + N$.

We use similar techniques to those in Section 4 to bound the PC angle $\Theta(\widehat{\boldsymbol{\Gamma}}_{\text{PC}}, \boldsymbol{\Gamma})$. We regard the term $N$ as a perturbation of order $\sigma^2$ of the fitted sample covariance matrix $\widehat{\mathbb{X}}^T \widehat{\mathbb{X}}$, and thus regard $\widehat{\boldsymbol{\Gamma}}_{\text{PC}}$ as a perturbation of $\widehat{\boldsymbol{\Gamma}}_{\text{PFC}}$.

**Definition 5.2.** *Writing $\lambda_i(\mathbf{X})$ for the ith ordered eigenvalue of $\mathbf{X}$, we define*

$$M = \min_{i \leq r} \left( \lambda_i(\widehat{\mathbb{X}}^T \widehat{\mathbb{X}}) - \lambda_{i+1}(\widehat{\mathbb{X}}^T \widehat{\mathbb{X}}) \right) \tag{9}$$

*for the minimum level spacing (this includes the spacing between zero and the non-zero eigenvalues, since $\lambda_{r+1}(\widehat{\mathbb{X}}^T \widehat{\mathbb{X}}) = 0$). We define the second level crossing event by*

$$\mathcal{L}_2 = \left\{ \lambda_1((\mathbb{X} - \widehat{\mathbb{X}})^T (\mathbb{X} - \widehat{\mathbb{X}})) \geq M \right\}.$$

**Lemma 5.3.** *If the errors $\boldsymbol{\epsilon}$ are independent standard normals, and if $\sqrt{M}/\sigma > \sqrt{n} + \sqrt{p}$ then the probability*

$$\mathbb{P}(\mathcal{L}_2) \leq \exp\left( -\frac{1}{2} \left( \sqrt{M}/\sigma - \sqrt{n} - \sqrt{p} \right)^2 \right).$$

*If $\mathcal{L}_2$ does not occur, then there are exactly r eigenvalues of $\mathbb{X}^T \mathbb{X}$ larger than $\lambda_r(\widehat{\mathbb{X}}^T \widehat{\mathbb{X}})$, so the r largest eigenvalues must correspond to the perturbed values of the original.*

We now prove bounds on the angle $\Theta(\widehat{\boldsymbol{\Gamma}}_{\text{PC}}, \widehat{\boldsymbol{\Gamma}}_{\text{PFC}})$ between the PC and PFC directions, where $C(\widehat{\boldsymbol{\Gamma}}_{\text{PC}}, \widehat{\boldsymbol{\Gamma}}_{\text{PFC}}) = \cot^2 \Theta(\widehat{\boldsymbol{\Gamma}}_{\text{PC}}, \widehat{\boldsymbol{\Gamma}}_{\text{PFC}})$ as before. Since $M$ and $\|(\mathbb{X} - \widehat{\mathbb{X}})^T (\mathbb{X} - \widehat{\mathbb{X}})\|$ are $O_{\mathbb{P}}(n)$, this angle is again of the order $O_{\mathbb{P}}(1/\sqrt{n})$.

**Proposition 5.4.** *In the case $d = r$, if the errors $\epsilon$ are independent standard normals, the conditional expectation*

$$\mathbb{E}\left(\Theta(\widehat{\mathbf{\Gamma}}_{\mathrm{PC}}, \widehat{\mathbf{\Gamma}}_{\mathrm{PFC}}) \Big| \mathcal{L}_2\right) \leq \arctan\left(\frac{\sigma^2\sqrt{np}}{M - \|(\mathbb{X} - \widehat{\mathbb{X}})^T(\mathbb{X} - \widehat{\mathbb{X}})\|}\right), \qquad (10)$$

*where, as in Equation (9), $M$ is the minimum eigenvalue spacing of $\widehat{\mathbb{X}}^T\widehat{\mathbb{X}}$.*

Using Proposition 5.4, we continue to explain the simulation graphs given in Section 5 of Cook [1], consider the means of the angles $\Theta(\widehat{\mathbf{\Gamma}}_{\mathrm{PC}}, \mathbf{\Gamma})$. Recall that Example 2.5, based on Figure 1 of [1], considers normally distributed errors $\epsilon$ with $p = 10$, $d = r = 1$, $\boldsymbol{\beta} = 1$, and varies parameters $n$, $\sigma$ and $\sigma_Y$. For simplicity, we replace some terms by their asymptotic limits to obtain more heuristic results. We recall that $\boldsymbol{\beta}\mathbf{F}^T\mathbf{F}\boldsymbol{\beta}^T \sim \sigma_Y^2 \chi_{n-1}^2$, so that asymptotically we can take $M = n\sigma_Y^2$. Similarly, we use the asymptotics given by Geman [6] and replace $\|(\mathbb{X} - \widehat{\mathbb{X}})^T(\mathbb{X} - \widehat{\mathbb{X}})\|$ by $\sigma^2 n$. Combining Lemma 5.3 and Proposition 5.4, we know that writing $C$ for $\arctan(\infty)$ (so that $C = \pi/2$ radians, or $90°$ for the graphs plotted) gives

$$\begin{aligned}
&\mathbb{E}\Theta(\widehat{\mathbf{\Gamma}}_{\mathrm{PC}}, \widehat{\mathbf{\Gamma}}_{\mathrm{PFC}}) \\
&= \mathbb{E}\left(\Theta(\widehat{\mathbf{\Gamma}}_{\mathrm{PC}}, \widehat{\mathbf{\Gamma}}_{\mathrm{PFC}}) \Big| \mathcal{L}_2\right)\mathbb{P}(\mathcal{L}_2) + \mathbb{E}\left(\Theta(\widehat{\mathbf{\Gamma}}_{\mathrm{PC}}, \widehat{\mathbf{\Gamma}}_{\mathrm{PFC}}) \Big| \mathcal{L}_2^c\right)\mathbb{P}(\mathcal{L}_2^c) \\
&\leq C\exp\left(-\frac{1}{2}\left(\frac{\sqrt{M}}{\sigma} - \sqrt{n} - \sqrt{p}\right)^2\right) \\
&\quad + \arctan\left(\frac{\sigma^2\sqrt{np}}{M - \|(\mathbb{X} - \widehat{\mathbb{X}})^T(\mathbb{X} - \widehat{\mathbb{X}})\|}\right) \\
&\leq C\exp\left(-\frac{1}{2}\left(\sqrt{n}(\sigma_Y/\sigma - 1) - \sqrt{p}\right)^2\right) + \arctan\left(\frac{\sigma^2\sqrt{p}}{\sqrt{n}(\sigma_Y^2 - \sigma^2)}\right).
\end{aligned}$$

We ignore the first term (since it decays exponentially fast in $n$ if $\sigma_Y^2 > \sigma^2$). In the case $r = d = 1$, angles satisfy a triangle inequality: that is, if the angle between vectors $\widehat{\mathbf{\Gamma}}_{\mathrm{PFC}}$ and $\mathbf{\Gamma}$ is $\theta_1$, and the angle between $\widehat{\mathbf{\Gamma}}_{\mathrm{PC}}$ and $\widehat{\mathbf{\Gamma}}_{\mathrm{PFC}}$ is $\theta_2$, then the angle between $\widehat{\mathbf{\Gamma}}_{\mathrm{PC}}$ and $\mathbf{\Gamma}$ is $\leq \theta_1 + \theta_2$. This means that, conditional on $\Theta(\widehat{\mathbf{\Gamma}}_{\mathrm{PFC}}, \mathbf{\Gamma})$, the angle between the PC directions and the true $\mathbf{\Gamma}$ is bounded above by

$$\mathbb{E}\Theta(\widehat{\mathbf{\Gamma}}_{\mathrm{PC}}, \mathbf{\Gamma}) \leq \Theta(\widehat{\mathbf{\Gamma}}_{\mathrm{PFC}}, \mathbf{\Gamma}) + \arctan\left(\frac{\sigma^2\sqrt{p}}{\sqrt{n}(\sigma_Y^2 - \sigma^2)}\right), \qquad (11)$$

where the distribution of $C = \cot^2\Theta(\widehat{\mathbf{\Gamma}}_{\mathrm{PFC}}, \mathbf{\Gamma})$ is given in Theorem 2.4. In each graph in Figure 2, we plot four curves, as follows:

1. We simulate directly from the distributions given by Theorem 2.4 and plot the sample mean of the PFC angle $\Theta(\widehat{\mathbf{\Gamma}}_{\mathrm{PFC}}, \mathbf{\Gamma})$ as ∘, as in Figure 1.

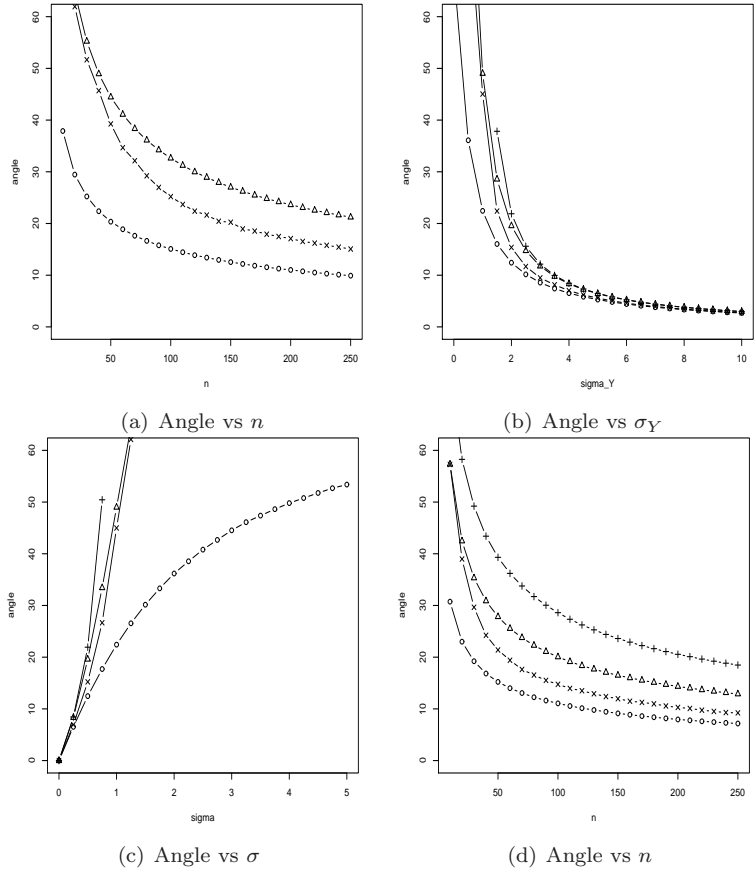

FIG 2. *Simulation of angle between $\mathbf{\Gamma}$ and both $\widehat{\mathbf{\Gamma}}_{\mathrm{PFC}}$ and $\widehat{\mathbf{\Gamma}}_{\mathrm{PC}}$, in the setting of Example 2.5. We plot the sample mean of the PFC angle $\Theta(\widehat{\mathbf{\Gamma}}_{\mathrm{PFC}}, \mathbf{\Gamma})$ as $\circ$, the true PC angle $\Theta(\widehat{\mathbf{\Gamma}}_{\mathrm{PC}}, \mathbf{\Gamma})$ as $\times$, the bound on the expected angle from Equation (11) as $+$ and the approximate expected angle from Equation (12) as $\triangle$. (a) Angle vs n, with $\sigma_Y = \sigma = 1$; (b) Angle vs $\sigma_Y$, with $n = 40$, $\sigma = 1$; (c) Angle vs $\sigma$, with $n = 40$, $\sigma_Y = 1$ (d) Angle vs n, with $\sigma_Y = \sqrt{2}$, $\sigma = 1$.*

2. Next, we plot the sample mean of the PC angle (based on direct simulation of 2500 samples from the underlying distribution) as $\times$.
3. We plot the bound given by the right hand side of (11) as $+$, noting that it is only useful when $\sigma_Y^2 > \sigma^2$.
4. Finally, in the spirit of Sibson [13], we use $\triangle$ to plot the angle corresponding to the leading term in $\sigma$ in the power series expansion in Equation (36) of the proof of Proposition 5.4, that is

$$\Theta(\widehat{\mathbf{\Gamma}}_{\mathrm{PFC}}, \mathbf{\Gamma}) + \arctan\left(\frac{\sigma^2 \sqrt{p}}{\sqrt{n} \sigma_Y^2}\right). \tag{12}$$

Providing a more complete description of the distribution of $\Theta(\widehat{\mathbf{\Gamma}}_{\mathrm{PC}}, \mathbf{\Gamma})$ than that given in Proposition 5.4, and extending such results to the case $d < r$, would require more advanced results from random matrix theory. However, the results shown in Figure 2 give a good explanation of Figure 1 of [1], and indicate parameter regions where the PFC and PC algorithms give close or differing results. In particular, the similarity of the plots in Figure 2 given by $\times$ and $\triangle$ shows that Equation (12) provides an accurate approximation to the PC performance. (In fact, (12) should approximately give an upper bound to the performance of this algorithm, because of the Jensen inequality used in Equation (37) of the proof of Proposition 5.4. However, it appears that the resulting distribution is sufficiently concentrated around its mean that the Jensen approximation is asymptotically accurate).

Since the bound given by (11) is only valid where $\sigma_Y^2 > \sigma^2$, we cannot plot it for all ranges of parameters considered by Cook – not at all in Figure 2(a), and only in a small region in Figure 2(c). We provide an extra series, Figure 2(d), where $\sigma_Y^2 = 2$, to show the dependence of this bound on $n$.

## Appendix A: Proofs

*Proof of Lemma 2.2.* As in [1], we write $\mathbf{1}_n$ for a $n \times 1$ matrix of 1s, so that if $\mathbf{X}$ is a $n \times p$ matrix with rows given by $\mathbf{X}_y^T$ then $\frac{1}{n}\mathbf{1}_n\mathbf{1}_n^T\mathbf{X}$ gives a matrix with rows all equal to the sample mean $\overline{\mathbf{X}}$. The assumption that $\sum_y \mathbf{f}_y = 0$ means that $\mathbf{1}_n^T\mathbf{F} = \mathbf{0}$, which implies that $\mathbf{P_F}\mathbf{1}_n = \mathbf{F}(\mathbf{F}^T\mathbf{F})^{-1}\mathbf{F}^T\mathbf{1}_n = \mathbf{0}$ – that is, $\frac{1}{n}\mathbf{1}_n\mathbf{1}_n^T$ and $\mathbf{P_F}$ represent projections onto orthogonal subspaces. Hence the matrix $\mathbb{X}$ defined to have rows $(\mathbf{X}_y - \overline{\mathbf{X}})^T$ can be expressed as $(\mathbf{I}_n - \frac{1}{n}\mathbf{1}_n\mathbf{1}_n^T)\mathbf{X}$ or, in terms of the quantities in Equation (1), as

$$\mathbb{X} = \left(\mathbf{I}_n - \frac{1}{n}\mathbf{1}_n\mathbf{1}_n^T\right)\left(\mathbf{1}_n\boldsymbol{\mu}^T + \mathbf{F}\boldsymbol{\beta}^T\boldsymbol{\Gamma}^T + \mathbf{E}\right) = \mathbf{F}\boldsymbol{\beta}^T\boldsymbol{\Gamma}^T + \left(\mathbf{I}_n - \frac{1}{n}\mathbf{1}_n\mathbf{1}_n^T\right)\mathbf{E}. \tag{13}$$

This means that the fitted matrix of predictors

$$\widehat{\mathbb{X}} = \mathbf{F}\boldsymbol{\beta}^T\boldsymbol{\Gamma}^T + \mathbf{P_F}\mathbf{E} = \mathbf{F}\left(\boldsymbol{\beta}^T\boldsymbol{\Gamma}^T + (\mathbf{F}^T\mathbf{F})^{-1}\mathbf{F}^T\mathbf{E}\right). \tag{14}$$

Using Equation (14), we rewrite $\widehat{\mathbb{X}} = \mathbf{FN}$ where $\mathbf{N} = (\boldsymbol{\beta}^T\boldsymbol{\Gamma}^T + (\mathbf{F}^T\mathbf{F})^{-1}\mathbf{F}^T\mathbf{E})$, so that we can define $\mathbf{V}^T = (\mathbf{F}^T\mathbf{F})^{-1/2}\mathbf{F}^T\widehat{\mathbb{X}} = (\mathbf{F}^T\mathbf{F})^{1/2}\mathbf{N}$, and

$$\widehat{\mathbb{X}}^T\widehat{\mathbb{X}} = \mathbf{N}^T\mathbf{F}^T\mathbf{F}\mathbf{N} = \left((\mathbf{F}^T\mathbf{F})^{1/2}\mathbf{N}\right)^T\left((\mathbf{F}^T\mathbf{F})^{1/2}\mathbf{N}\right) = \mathbf{V}\mathbf{V}^T. \tag{15}$$

Any vector $\mathbf{w}$ orthogonal to the span of the columns of $\mathbf{V}$ is a 0-eigenvector of $\widehat{\mathbb{X}}^T\widehat{\mathbb{X}}$, since for such a vector Equation (15) shows that

$$\widehat{\mathbb{X}}^T\widehat{\mathbb{X}}\mathbf{w} = \mathbf{V}\mathbf{V}^T\mathbf{w} = \mathbf{V}\mathbf{0} = \mathbf{0},$$

and similarly the span of the columns of $\mathbf{V}$ is preserved by $\widehat{\mathbb{X}}^T\widehat{\mathbb{X}}$. Hence the $r$ columns of $\mathbf{V}$ span the same space as the eigenvectors of $\widehat{\mathbb{X}}^T\widehat{\mathbb{X}}$ with the $r$ largest eigenvalues. □

*Proof of Lemma 2.3.* For any matrix $\mathbf{A}$, the entries of $\mathbf{AE}$ have mean $\mathbb{E}(\mathbf{AE}) = \mathbf{A}\mathbb{E}(\mathbf{E}) = \sigma \mathbf{A}\mathbb{E}(\boldsymbol{\epsilon}) = 0$, allowing us to deduce the mean of $\mathbf{V}$. Similarly, we use the well-known fact that for any $\mathbf{A}$ and $\mathbf{B}$:

$$\text{Cov}\,((\mathbf{AE})_{ij}, (\mathbf{BE})_{kl}) = \sum_{r,s} \mathbf{A}_{ir}\mathbf{B}_{ks}\text{Cov}\,(\mathbf{E}_{rj}, \mathbf{E}_{sl}) = \sigma^2(\mathbf{AB}^T)_{ik}\delta_{jl}. \quad (16)$$

Taking $\mathbf{A} = \mathbf{B} = (\mathbf{F}^T\mathbf{F})^{-1/2}\mathbf{F}^T$ means that $\mathbf{AB}^T = \mathbf{I}$, so by Equation (16) the entries of $\mathbf{V}$ satisfy

$$\text{Cov}\,(\mathbf{V}_{ij}, \mathbf{V}_{kl}) = \text{Cov}\,((\mathbf{F}^T\mathbf{F})^{-1/2}\mathbf{F}^T\mathbf{E})_{ij}, (\mathbf{F}^T\mathbf{F})^{-1/2}\mathbf{F}^T\mathbf{E})_{kl}) = \sigma^2\delta_{ik}\delta_{jl},$$

and the result follows. □

*Proof of Theorem 2.4.* In the case $d = r$, Lemma 2.2 tells us that we can choose to take a PFC estimate $\widehat{\boldsymbol{\Gamma}}_{\text{PFC}}$ given by an orthogonal transformation of the columns of $\mathbf{V}$ defined in Equation (3). Write $\mathbf{b}^{(i)}$, where $i = 1, \ldots, d$, for the columns of $\boldsymbol{\Gamma}$, which form an orthonormal set since, by assumption, $\boldsymbol{\Gamma}^T\boldsymbol{\Gamma} = \mathbf{I}_d$. We extend this to create an orthonormal basis $\{\mathbf{b}^{(1)}, \ldots, \mathbf{b}^{(p)}\}$ for the whole of $\mathbb{R}^p$, and write $\mathbf{G}$ for the $p \times p$ matrix made up of the complete set of columns $\mathbf{b}^{(i)}$, with $\mathbf{G}^T\mathbf{G} = \mathbf{I}_p$.

We express $\mathbf{V}$ in this new basis, for $k = 1, \ldots, r$ we expand the $k$th column of $\mathbf{V}$ as:

$$\mathbf{V}^{(k)} = \sum_{i=1}^{p} \mathbf{A}_{ik}\mathbf{b}^{(i)}, \quad \text{where } \mathbf{A}_{ik} = (\mathbf{b}^{(i)})^T\mathbf{V}^{(k)}. \quad (17)$$

Equivalently, we write the $p \times r$ matrix $\mathbf{A} = \mathbf{G}^T\mathbf{V} = \boldsymbol{\mu} + \mathbf{C}$. Here $\boldsymbol{\mu} = \mathbf{G}^T\boldsymbol{\Gamma}\boldsymbol{\beta}(\mathbf{F}^T\mathbf{F})^{1/2}$, consists of a $d \times r$ block of the form $\boldsymbol{\beta}(\mathbf{F}^T\mathbf{F})^{1/2}$ and a $(p-d) \times r$ zero block, and $\mathbf{C} = \mathbf{G}^T\mathbf{E}^T\mathbf{F}(\mathbf{F}^T\mathbf{F})^{-1/2}$. Equation (17) and Lemma 2.3 show that $\mathbf{A}$ has mean $\boldsymbol{\mu}$, so that for $i \geq d+1$ the $\mathbf{A}_{ik}$ has mean 0, whereas for $i \leq d$ the $\mathbf{A}_{ik}$ has mean

$$\boldsymbol{\mu}_{ik} = (\boldsymbol{\beta}(\mathbf{F}^T\mathbf{F})^{1/2})_{ik}. \quad (18)$$

Using Lemma 2.3, we know that for any $1 \leq i, j \leq p$ and $1 \leq k, l \leq r$:

$$\begin{aligned}
\text{Cov}\,(\mathbf{A}_{ik}, \mathbf{A}_{jl}) &= \sum_{u=1}^{p}\sum_{v=1}^{p}(\mathbf{b}^{(i)})_u(\mathbf{b}^{(j)})_v\text{Cov}\,(\mathbf{V}_{uk}, \mathbf{V}_{vl}) \\
&= \sigma^2\sum_{u=1}^{p}\sum_{v=1}^{p}(\mathbf{b}^{(i)})_u(\mathbf{b}^{(j)})_v\delta_{uv}\delta_{kl} = \sigma^2\delta_{kl}\delta_{ij}. \quad (19)
\end{aligned}$$

The fact that $\boldsymbol{\Gamma}^T\boldsymbol{\Gamma} = \mathbf{I}_d$ implies that the projection $\mathbf{P}_{\boldsymbol{\Gamma}} = \boldsymbol{\Gamma}\boldsymbol{\Gamma}^T$, so that the projections of the $k$th column become

$$\begin{aligned}
\mathbf{P}_{\boldsymbol{\Gamma}}\mathbf{V}^{(k)} &= \boldsymbol{\Gamma}\boldsymbol{\Gamma}^T\mathbf{V}^{(k)} = \sum_{i=1}^{d}\mathbf{A}_{ik}\mathbf{b}^{(i)}, \\
(\mathbf{I} - \mathbf{P}_{\boldsymbol{\Gamma}})\mathbf{V}^{(k)} &= (\mathbf{I}_p - \boldsymbol{\Gamma}\boldsymbol{\Gamma}^T)\mathbf{V}^{(k)} = \sum_{i=d+1}^{p}\mathbf{A}_{ik}\mathbf{b}^{(i)}.
\end{aligned}$$

Hence, the respective Frobenius norms become

$$\|\mathbf{P}_{\boldsymbol{\Gamma}}\mathbf{V}\|_F^2 = \sum_{k=1}^{r}\left|\sum_{i=1}^{d}\mathbf{A}_{ik}\mathbf{b}^{(i)}\right|^2 = \sum_{k=1}^{r}\sum_{i=1}^{d}\mathbf{A}_{ik}^2, \qquad (20)$$

$$\|(\mathbf{I}_p - \mathbf{P}_{\boldsymbol{\Gamma}})\mathbf{V}\|_F^2 = \sum_{k=1}^{r}\left|\sum_{i=d+1}^{p}\mathbf{A}_{ik}\mathbf{b}^{(i)}\right|^2 = \sum_{k=1}^{r}\sum_{i=d+1}^{p}\mathbf{A}_{ik}^2. \qquad (21)$$

Note that the arguments so far do not require the assumption that the errors $\boldsymbol{\epsilon}$ are normal. Under this additional assumption, we deduce that the $\mathbf{A}_{ik}$ are normal and independent, with common variance $\sigma^2$. In particular, under this assumption the expressions (20) and (21) are independent of each other.

Further, if the errors are normal, Equation (21) means that $\|(\mathbf{I}-\mathbf{P}_{\boldsymbol{\Gamma}})\mathbf{V}\|_F^2/\sigma^2$ has a central $\chi^2$ distribution with $r(p-d)$ degrees of freedom. Similarly, Equation (20) means that $\|\mathbf{P}_{\boldsymbol{\Gamma}}\mathbf{V}\|_F^2/\sigma^2$ has a non-central $\chi^2$ distribution with $rd$ degrees of freedom, and non-centrality parameter

$$\Lambda = \frac{1}{\sigma^2}\sum_{k=1}^{r}\sum_{i=1}^{d}\boldsymbol{\mu}_{ik}^2 = \frac{1}{\sigma^2}\text{tr}\,(\boldsymbol{\beta}(\mathbf{F}^T\mathbf{F})\boldsymbol{\beta}^T), \qquad (22)$$

and the proof is complete. □

*Proof of Lemma 3.1.* We write $\mathbf{C}_{ik} = \mathbf{A}_{ik} - \boldsymbol{\mu}_{ik}$, and use Equations (18)–(21). In particular, $\mathbf{C}_{ik}$ have mean zero and variance $\sigma^2$, and Equation (18) implies that $\sum_{k,i}\boldsymbol{\mu}_{ik}^2 = \text{tr}\,(\boldsymbol{\beta}(\mathbf{F}^T\mathbf{F})\boldsymbol{\beta}^T)$. Hence we know that $N = \sum_{k=1}^{r}\sum_{i=1}^{d}\mathbf{A}_{ik}^2$ has mean $\text{tr}\,(\boldsymbol{\beta}(\mathbf{F}^T\mathbf{F})\boldsymbol{\beta}^T) + rd\sigma^2$ and $D = \sum_{k=1}^{r}\sum_{i=d+1}^{p}\mathbf{A}_{ik}^2$ has mean $r(p-d)\sigma^2$.

For $i \leq d$, we write $\mathbf{C}_{ik} = \sum_{r,s}\boldsymbol{\Gamma}_{ri}\mathbf{E}_{sr}\mathbf{W}_{sk}$ where $\mathbf{W} = \mathbf{F}(\mathbf{F}^T\mathbf{F})^{-1/2}$, and $\mathbf{M} = \boldsymbol{\Gamma}\boldsymbol{\Gamma}^T$. We deduce that $\text{Var}\,(N)$ equals

$$\sum_{k=1}^{r}\left(\mathbb{E}\left(\sum_{i=1}^{d}\mathbf{C}_{ik}^2\right)^2 + 4\sum_{i,j=1}^{d}\boldsymbol{\mu}_{ik}\mathbb{E}\mathbf{C}_{ik}\mathbf{C}_{kj}^2 - d^2\sigma^4\right) + 4\text{tr}\,(\boldsymbol{\beta}(\mathbf{F}^T\mathbf{F})\boldsymbol{\beta}^T)\sigma^2. \quad (23)$$

In this expansion, we know that $\mathbb{E}\mathbf{C}_{ik}\mathbf{C}_{jk}^2 = 0$, since it can be written as a linear combination of expectations of terms in $\mathbf{E}_{ij}$ of degree 3. Each such term contains a term of odd degree, so independence and symmetry means the resulting expectation is zero. Similarly, using independence and symmetry, since only terms of the form $\mathbb{E}\mathbf{E}_{ij}^4$ and $\mathbb{E}\mathbf{E}_{ij}^2\mathbf{E}_{kl}^2$ make a non-zero contribution, for any

$k$ we can expand

$$\mathbb{E}\left(\sum_{i=1}^d \mathbf{C}_{ik}^2\right)^2$$
$$= \mathbb{E}\left(\sum_{r,s,t,u} \mathbf{W}_{sk}\mathbf{W}_{tk}\mathbf{E}_{sr}\mathbf{E}_{tu}\mathbf{M}_{ru}\right)^2$$
$$= \sigma^4 \sum_{r,s} \mathbf{W}_{sk}^4 \mathbf{M}_{rr}^2(m_4-3) + \sigma^4 \sum_{r,s,t,u} \mathbf{W}_{sk}^2\mathbf{W}_{tk}^2(\mathbf{M}_{rr}\mathbf{M}_{uu} + 2\mathbf{M}_{ru}\mathbf{M}_{ur})$$
$$\leq \sigma^4 \left(\sum_{r,s,t,u} \mathbf{W}_{sk}^2\mathbf{W}_{tk}^2 \mathbf{M}_{rr}\mathbf{M}_{uu} + (m_4-1)\mathbf{M}_{ru}\mathbf{M}_{ur}\right)$$
$$= \sigma^4\left(d^2 + (m_4-1)d\right)$$

since $\sum_s \mathbf{W}_{sk}^2 = \sum_s \mathbf{W}_{ks}^T \mathbf{W}_{sk} = 1$, the $\sum_r \mathbf{M}_{rr} = \operatorname{tr}(\mathbf{\Gamma}\mathbf{\Gamma}^T) = \operatorname{tr}(\mathbf{\Gamma}^T\mathbf{\Gamma}) = d$, and $\sum_u \mathbf{M}_{ru}\mathbf{M}_{ur} = \mathbf{M}_{rr}$, so the result follows. Similarly we deduce $\operatorname{Var}(D) \leq (p-d)T$. □

If $\boldsymbol{\epsilon}$ are independent with mean zero and finite 4th moment, but no longer symmetric, we can prove similar results. In this case, the $\sum_{i,j=1}^d \boldsymbol{\mu}_{ik}\mathbb{E}\mathbf{C}_{ik}\mathbf{C}_{jk}^2$ term becomes $O_{\mathbb{P}}(\sqrt{n})$, which allows this proof to be adapted.

Fulton [5] reviews many facts concerning the eigenvalues of sums of matrices, including the following result:

**Lemma A.1** (Weyl [14]). *For any real symmetric* $\mathbf{B}$ *and* $\mathbf{L}$:

$$\lambda_i(\mathbf{B}) + \lambda_p(\mathbf{L}) \leq \lambda_i(\mathbf{B}+\mathbf{L}) \leq \lambda_i(\mathbf{B}) + \lambda_1(\mathbf{L}) \quad \text{for } i=1,2,\ldots,p.$$

*Proof of Lemma 3.2.* Since $\mathbf{F}$ has rows which correspond to independent samples of $\mathbf{Y}$, the Law of Large Numbers means that $(\mathbf{F}^T\mathbf{F})/n$ converges elementwise almost surely to a real symmetric matrix $\mathbf{\Phi}$. Hence, using the union bound, $\sup_{i,j}\left|(\boldsymbol{\beta}(\mathbf{F}^T\mathbf{F})\boldsymbol{\beta}^T)_{ij}/n - \mathbf{\Phi}_{ij}\right|$ can be made arbitrarily small for all $n$ sufficiently large with probability $1-\epsilon$. This implies that the spectral norm of the difference $\|(\boldsymbol{\beta}(\mathbf{F}^T\mathbf{F})\boldsymbol{\beta}^T/n - \mathbf{\Phi}\|$ can be made smaller than $\delta$ with the same probability. By Lemma A.1, this implies that the $i$th eigenvalue $\lambda_i(\boldsymbol{\beta}(\mathbf{F}^T\mathbf{F})\boldsymbol{\beta}^T/n)$ is arbitrarily close to $\phi_i = \lambda_i(\mathbf{\Phi})$. □

*Proof of Theorem 3.3.* We write $C(\widehat{\mathbf{\Gamma}}_{\text{PFC}}, \mathbf{\Gamma}) = N/D$, in the notation of Lemma 3.1 and define the event

$$B_{K_1,K_2,n^*} = \left\{K_1 \leq \frac{\operatorname{tr}(\boldsymbol{\beta}(\mathbf{F}^T\mathbf{F})\boldsymbol{\beta}^T)}{n} \leq K_2 \text{ for all } n \geq n^*\right\}.$$

By Lemma 3.2, given any $\alpha$ there exist $K_1, K_2, n^*$ such that $\mathbb{P}(B_{K_1,K_2,n^*}) \geq 1 - \alpha/3$. Note that on $B_{K_1,K_2,n^*}$, Lemma 3.1 implies that the mean $\mathbb{E}N \geq$ $K_1 n$ and Var $N \leq Td + 4K_2 n$. We use Chebyshev and choose $N_+^* = K_1 n - \sqrt{3(Td + 4K_2 n)/\alpha}$ to write that

$$\mathbb{P}(N \leq N_+^* | B_{K_1, K_2, n^*}) \leq \frac{Td + 4K_2 n}{(N^* - K_1 n)^2} \leq \frac{\alpha}{3}. \quad (24)$$

We bound an upper confidence interval by taking $X_+ = (3r(p-d)\sigma^2/\alpha)$ and combining the fact that $D$ is independent of the event $B_{K_1, K_2, n^*}$ with Markov's inequality to obtain

$$\begin{aligned}
&\mathbb{P}(\Theta(\widehat{\mathbf{\Gamma}}_{\text{PFC}}, \mathbf{\Gamma}) \geq \Theta_+^*) \\
&= \mathbb{P}(\Theta(\widehat{\mathbf{\Gamma}}_{\text{PFC}}, \mathbf{\Gamma}) \geq \Theta_+^* | B_{K_1, K_2, n^*}) \mathbb{P}(B_{K_1, K_2, n^*}) \\
&\quad + \mathbb{P}(\Theta(\widehat{\mathbf{\Gamma}}_{\text{PFC}}, \mathbf{\Gamma}) \geq \Theta_+^* | B_{K_1, K_2, n^*}^c) \mathbb{P}(B_{K_1, K_2, n^*}^c) \\
&\leq \mathbb{P}\left(D/N \geq \tan^2(\Theta_+^*) | B_{K_1, K_2, n^*}\right) + \mathbb{P}(B_{K_1, K_2, n^*}^c) \\
&\leq \mathbb{P}\left(D/N \geq \tan^2(\Theta_+^*), D \leq X_+ | B_{K_1, K_2, n^*}\right) \\
&\quad + \mathbb{P}\left(D/N \geq \tan^2(\Theta_+^*), D \geq X_+ | B_{K_1, K_2, n^*}\right) + \alpha/3 \\
&\leq \mathbb{P}\left(N \leq X_+/\tan^2(\Theta_+^*) | B_{K_1, K_2, n^*}\right) + \mathbb{P}\left(D \geq X_+ | B_{K_1, K_2, n^*}\right) + \alpha/3 \\
&\leq \mathbb{P}\left(N \leq X_+/\tan^2(\Theta_+^*)\right) + 2\alpha/3. \quad (25)
\end{aligned}$$

We obtain the required bounds by equating $N_+^*$ and $X_+/\tan^2(\Theta_+^*)$, and using Equation (24) and the fact that $\arctan(t) \leq t$. That is, we can take

$$\Theta_+^*(\alpha) = \arctan\left(\sqrt{\frac{X_+}{K_1 n - \sqrt{3(T + 4K_2 n)/\alpha}}}\right), \quad (26)$$

and substitute the value $X_+ = (3r(p-d)\sigma^2/\alpha)$ given above. We can find a lower confidence interval, since for any $X_-$, a similar argument to Equation (25) gives

$$\mathbb{P}(\Theta(\widehat{\mathbf{\Gamma}}_{\text{PFC}}, \mathbf{\Gamma}) \leq \Theta_-^*) \leq \mathbb{P}(D \leq X_-) + \mathbb{P}\left(N \geq X_-/\tan^2(\Theta_-^*)\right) + \mathbb{P}(B_{K_1, K_2, n^*}^c).$$

By Chebyshev and Lemma 3.1, we can take $X_- = \sigma^2 r(p-d) - \sqrt{3T(p-d)/\alpha}$ to ensure that $\mathbb{P}(D \leq X_-) \leq \alpha/3$, and $N_-^* = K_2 n + r d \sigma^2 + \sqrt{3(Td + 4K_2 n)/\alpha}$ to ensure that $\mathbb{P}(N \geq N_-^*) \leq \alpha/3$. Again, equating $N_-^*$ and $X_-/\tan^2(\Theta_-^*)$, we deduce that

$$\Theta_-^*(\alpha) = \arctan\left(\sqrt{\frac{X_-}{K_2 n + r d \sigma^2 + \sqrt{3(T + 4K_2 n)/\alpha}}}\right), \quad (27)$$

and the result follows in the same way. □

*Proof of Theorem 4.3.* During this proof, for simplicity, we write $\lambda_i$ for $\lambda_i(\mathbf{U}\mathbf{U}^T)$, which equals $\lambda_i(\boldsymbol{\beta}(\mathbf{F}^T \mathbf{F})\boldsymbol{\beta}^T)$ if $i \leq d$ and zero otherwise. First, we introduce the event

$$C_{\mathbf{\Phi}, n^*} = \left\{\lambda_i(\mathbf{\Phi}) - \delta \leq \frac{\lambda_i}{n} \leq \lambda_i(\mathbf{\Phi}) + \delta \text{ for all } i = 1, \ldots, d, \text{ and all } n \geq n^*\right\},$$

where $\delta$ is defined in terms of the interlevel spacing as $\delta = \min_{i \leq d}(\lambda_i(\mathbf{\Phi}) - \lambda_{i+1}(\mathbf{\Phi}))/10$, writing $\lambda_{d+1}(\mathbf{\Phi}) = 0$. Note that $\delta > 0$ by assumption.

We mirror the proof of Theorem 3.3 by conditioning on whether $C_{\mathbf{\Phi},n^*}$ and $\mathcal{L}_1$ occur, to obtain

$$\begin{aligned}
&\mathbb{P}(\Theta(\widehat{\mathbf{\Gamma}}_{\text{PFC}}, \mathbf{\Gamma}) \geq \Theta^*) \\
&= \mathbb{P}(\Theta(\widehat{\mathbf{\Gamma}}_{\text{PFC}}, \mathbf{\Gamma}) \geq \Theta^* | C^c_{\mathbf{\Phi},n^*})\mathbb{P}(C^c_{\mathbf{\Phi},n^*}) \\
&\quad + \mathbb{P}(\Theta(\widehat{\mathbf{\Gamma}}_{\text{PFC}}, \mathbf{\Gamma}) \geq \Theta^* | \mathcal{L}_1, C_{\mathbf{\Phi},n^*})\mathbb{P}(\mathcal{L}_1|C_{\mathbf{\Phi},n^*})\mathbb{P}(C_{\mathbf{\Phi},n^*}) \\
&\quad + \mathbb{P}(\Theta(\widehat{\mathbf{\Gamma}}_{\text{PFC}}, \mathbf{\Gamma}) \geq \Theta^* | \mathcal{L}_1^c, C_{\mathbf{\Phi},n^*})\mathbb{P}(\mathcal{L}_1^c, C_{\mathbf{\Phi},n^*}) \\
&\leq \mathbb{P}(C^c_{\mathbf{\Phi},n^*}) + \mathbb{P}(\mathcal{L}_1|C_{\mathbf{\Phi},n^*}) + \mathbb{P}(\sin\Theta(\widehat{\mathbf{\Gamma}}_{\text{PFC}}, \mathbf{\Gamma}) \geq \sin\Theta^* | \mathcal{L}_1^c, C_{\mathbf{\Phi},n^*}). \quad (28)
\end{aligned}$$

Lemma 3.2 means that the first term in Equation (28) is less than $\alpha/3$ for $n^*$ sufficiently large. Next we bound the second term in Equation (28), using the fact that $\lambda_1(\mathbf{L}) - \lambda_p(\mathbf{L}) \leq 4\sigma\sqrt{\|\mathbf{SS}^T\|\|\mathbf{UU}^T\|} + \sigma^2\|\mathbf{SS}^T\|$, where we write $\|\cdot\|$ for the spectral norm. Completing the square and using Equation (8), we deduce that

$$\begin{aligned}
\mathbb{P}(\mathcal{L}_1|C_{\mathbf{\Phi},n^*}) &= \mathbb{P}(\lambda_1(\mathbf{L}) - \lambda_p(\mathbf{L}) \geq \lambda_d | C_{\mathbf{\Phi},n^*}) \\
&\leq \mathbb{P}\left(4\sigma\sqrt{\|\mathbf{SS}^T\|}\sqrt{\|\mathbf{UU}^T\|} + \sigma^2\|\mathbf{SS}^T\| \geq \lambda_d \bigg| C_{\mathbf{\Phi},n^*}\right) \\
&= \mathbb{P}\left(\sigma\sqrt{\|\mathbf{SS}^T\|} \geq \sqrt{4\lambda_1 + \lambda_d} - 2\sqrt{\lambda_1} \bigg| C_{\mathbf{\Phi},n^*}\right) \\
&= \mathbb{P}\left(\sqrt{\|\mathbf{SS}^T\|} \geq \sqrt{r} + \sqrt{p} + u \bigg| C_{\mathbf{\Phi},n^*}\right) \\
&\leq \exp(-u^2/2)
\end{aligned}$$

where $u = (\sqrt{4\lambda_1 + \lambda_d} - 2\sqrt{\lambda_1})/\sigma - \sqrt{r} - \sqrt{p}$. The last inequality uses the fact that $\mathbf{S}$ is a $p \times r$ matrix of independent standard Gaussians, so we can apply Equation (8). Conditioning on the event $C_{\mathbf{\Phi},n^*}$ ensures that

$$u \geq \sqrt{n}(\sqrt{4\lambda_1(\mathbf{\Phi}) + \lambda_d(\mathbf{\Phi}) - 5\delta} - 2\sqrt{\lambda_1(\mathbf{\Phi}) + \delta})/\sigma - \sqrt{r} - \sqrt{p} = a\sqrt{n} - \sqrt{r} - \sqrt{p},$$

where the choice of $\delta$ ensures that $a > 0$, so the second term in Equation (28) is less than $\alpha/3$ for $n$ sufficiently large.

Finally we condition on the event $\mathcal{L}_1^c \cap C_{\mathbf{\Phi},n^*}$ The key result is (as in [13]) that if $\mathbf{e}_i$ is a $\lambda_i$-eigenvector of $\mathbf{UU}^T$ then if $\mathbf{f}_i$ satisfies

$$\mathbf{f}_i = -(\mathbf{UU}^T - \lambda_i\mathbf{I}_p)^+(\mathbf{L} - \mu_i\mathbf{I}_p)(\mathbf{e}_i + \mathbf{f}_i), \quad (29)$$

then $\mathbf{e}_i + \mathbf{f}_i$ is a $(\lambda_i + \mu_i)$-eigenvector of $\widehat{\mathbb{X}}^T\widehat{\mathbb{X}} = \mathbf{UU}^T + \mathbf{L}$, with $\mathbf{f}_i \perp \mathbf{e}_i$. Hence, Equation (29) tells us that, conditioned on $\mathcal{L}_1^c$, these perturbed eigenvectors form the PFC directions. Since we condition on the event $C_{\mathbf{\Phi},n^*}$, the norm

$$\|(\mathbf{UU}^T - \lambda_i\mathbf{I}_p)^+\| \leq \frac{1}{\min_{j \neq i}(\lambda_j - \lambda_i)} \leq \frac{10}{n\delta} \quad \text{for all } i. \quad (30)$$

We need to probabilistically bound the norm of $(\mathbf{L} - \mu_i \mathbf{I}_p)$. Lemma A.1 implies that $\lambda_p(\mathbf{L}) \leq \mu_i \leq \lambda_1(\mathbf{L})$, so that for all $i$ the norm

$$\|(\mathbf{L} - \mu_i \mathbf{I}_p)\| \leq \lambda_1(\mathbf{L}) - \lambda_p(\mathbf{L}) \leq 4\sigma \sqrt{\|\mathbf{SS}^T\| \|\mathbf{UU}^T\|} + \sigma^2 \|\mathbf{SS}^T\|. \quad (31)$$

Using Proposition 4.1, we know that since $\mathbf{S}$ is a $p \times r$ matrix of standard Gaussians, there exists $K_1$ such that $\mathbb{P}(\|\mathbf{SS}^T\| \geq K_1) \leq \alpha/3$. Since we condition on the event $C_{\mathbf{\Phi}, n^*}$, the $\|\mathbf{UU}^T\| \leq (\lambda_1(\mathbf{\Phi}) + \delta)n = K_2 n$ for $n$ sufficiently large. Overall, we deduce that

$$\mathbb{P}\left(\|\mathbf{L} - \mu_i \mathbf{I}_p)\| \geq 4\sigma\sqrt{K_1 K_2 n} + \sigma^2 K_1 \,\Big|\, \mathcal{L}_1^c, C_{\mathbf{\Phi}, n^*}\right) \leq \alpha/3. \quad (32)$$

Taking $K = 10(4\sigma\sqrt{K_1 K_2} + \sigma^2 K_1)/\delta$ and using Equations (29) and (30), we deduce that, since $\sin\Theta(\widehat{\mathbf{\Gamma}}_{\text{PFC}}, \mathbf{\Gamma})$ represents a weighted average of sines of angles between PFC and unperturbed eigenvectors, it is dominated by the maximum sine of such an angle. In other words

$$\mathbb{P}(\sin\Theta(\widehat{\mathbf{\Gamma}}_{\text{PFC}}, \mathbf{\Gamma}) \geq K/\sqrt{n}\,\Big|\, \mathcal{L}_1^c, C_{\mathbf{\Phi}, n^*})$$
$$\leq \mathbb{P}\left(\max_i \left(\|(\mathbf{UU}^T - \lambda_i \mathbf{I}_p)^+ (\mathbf{L} - \mu_i \mathbf{I}_p)\|\right) \geq K/\sqrt{n}\,\Big|\, \mathcal{L}_1^c, C_{\mathbf{\Phi}, n^*}\right)$$
$$\leq \mathbb{P}\left(\max_i \|\mathbf{L} - \mu_i \mathbf{I}_p\| \geq \frac{\sqrt{n}\delta K}{10}\,\Big|\, \mathcal{L}_1^c, C_{\mathbf{\Phi}, n^*}\right)$$
$$= \mathbb{P}\left(\|\mathbf{L} - \mu_i \mathbf{I}_p)\| \geq (4\sigma\sqrt{K_1 K_2} + \sigma^2 K_1)\sqrt{n}\,\Big|\, \mathcal{L}_1^c, C_{\mathbf{\Phi}, n^*}\right)$$
$$\leq \mathbb{P}\left(\|\mathbf{L} - \mu_i \mathbf{I}_p)\| \geq 4\sigma\sqrt{K_1 K_2 n} + \sigma^2 K_1 \,\Big|\, \mathcal{L}_1^c, C_{\mathbf{\Phi}, n^*}\right) \leq \alpha/3, \quad (33)$$

by Equation (32). Taking $\sin\Theta^* = K/\sqrt{n}$, Equation (33) tells us that the third term of Equation (28) is less than $\alpha/3$. The result follows using the fact that $\Theta^* = \arcsin(K/\sqrt{n}) = O(1/\sqrt{n})$. □

*Proof of Lemma 5.1.* The key is to notice that using Equations (13) and (14), and writing $\mathbf{P_G} = \left(\mathbf{I}_n - \frac{1}{n}\mathbf{1}_n \mathbf{1}_n^T - \mathbf{P_F}\right)$ for the projection onto the space orthogonal to the columns of $\mathbf{F}$ and to the column made up of 1s, we can express

$$\mathbb{X} - \widehat{\mathbb{X}} = \left(\mathbf{I}_n - \frac{1}{n}\mathbf{1}_n \mathbf{1}_n^T - \mathbf{P_F}\right) \mathbf{E}, \quad (34)$$

which is not a function of $\mathbf{\Gamma}$. As before, since $\mathbf{P_F}$ and $\mathbf{P_G}$ represent projection onto orthogonal spaces, the $\mathbf{P_F P_G^T} = \mathbf{0}$. This means that

$$\widehat{\mathbb{X}}^T(\mathbb{X} - \widehat{\mathbb{X}}) = \mathbb{X}^T \mathbf{P_F P_G^T} \mathbb{X} = \mathbf{0},$$

so, since the cross-term vanishes,

$$\mathbb{X}^T \mathbb{X} = (\widehat{\mathbb{X}} + (\mathbb{X} - \widehat{\mathbb{X}}))^T (\widehat{\mathbb{X}} + (\mathbb{X} - \widehat{\mathbb{X}})) = \widehat{\mathbb{X}}^T \widehat{\mathbb{X}} + (\mathbb{X} - \widehat{\mathbb{X}})^T(\mathbb{X} - \widehat{\mathbb{X}}). \quad (35)$$

Further, by Equation (16), the relation $\mathbf{P_F P_G^T} = \mathbf{0}$ means that the terms $\widehat{\mathbb{X}}$ and $\mathbb{X} - \widehat{\mathbb{X}}$ are uncorrelated. Hence, if the errors $\boldsymbol{\epsilon}$ are normally distributed, the terms $\widehat{\mathbb{X}}$ and $\mathbb{X} - \widehat{\mathbb{X}}$ are independent and the lemma follows. □

*Proof of Lemma 5.3.* Note that Equation (34) means that $(\mathbb{X} - \widehat{\mathbb{X}})^T(\mathbb{X} - \widehat{\mathbb{X}}) = \sigma^2 \mathbf{X}^T \mathbf{P_G} \mathbf{X}$, where $\mathbf{X}$ is an $n \times p$ Wishart matrix with standard Gaussian entries independent of $\widehat{\mathbb{X}}^T \widehat{\mathbb{X}}$. We take $t = \sqrt{M}/\sigma - (\sqrt{n} + \sqrt{p})$ and apply the result of Davidson and Szarek, Equation (8).

Since $(\mathbb{X} - \widehat{\mathbb{X}})^T(\mathbb{X} - \widehat{\mathbb{X}})$ is positive-definite, $\lambda_p((\mathbb{X} - \widehat{\mathbb{X}})^T(\mathbb{X} - \widehat{\mathbb{X}})) \geq 0$, and since $\lambda_i(\widehat{\mathbb{X}}^T\widehat{\mathbb{X}}) = 0$ for $i \geq r+1$, Lemma A.1 and Equation (35) imply that

$$\lambda_i(\widehat{\mathbb{X}}^T\widehat{\mathbb{X}}) \leq \lambda_i(\mathbb{X}^T\mathbb{X}) \leq \lambda_i(\widehat{\mathbb{X}}^T\widehat{\mathbb{X}}) + \lambda_1((\mathbb{X} - \widehat{\mathbb{X}})^T(\mathbb{X} - \widehat{\mathbb{X}}))$$

Hence, if $\mathcal{L}_2$ does not occur, that is if $\lambda_1((\mathbb{X}-\widehat{\mathbb{X}})^T(\mathbb{X}-\widehat{\mathbb{X}})) < M$, there are exactly $r$ eigenvalues of $\mathbb{X}^T\mathbb{X}$ larger than $\lambda_r(\widehat{\mathbb{X}}^T\widehat{\mathbb{X}})$, and the $r$ largest eigenvalues must correspond to the perturbed values of the original. □

Next we prove another technical lemma:

**Lemma A.2.** *Given $n \times n$ projection matrix $\mathbf{P_G}$ and $n \times p$ matrix $\mathbf{X}$ with independent standard Gaussian entries, consider matrices $\mathbf{R} = \mathbf{X}\mathbf{W}$ and $\mathbf{S} = \mathbf{X}\mathbf{V}$, where $\mathbf{W}^T\mathbf{V} = \mathbf{0}$. Then*

$$\mathbb{E}\left(\mathrm{tr}\,(\mathbf{R}\mathbf{R}^T\mathbf{P_G}\mathbf{S}\mathbf{S}^T\mathbf{P_G})\right) = \mathrm{tr}\,(\mathbf{W}^T\mathbf{W})\mathrm{tr}\,(\mathbf{V}^T\mathbf{V})\mathrm{rank}\,(\mathbf{G}).$$

*Proof.* Using the representation $\mathbf{R} = \mathbf{X}\mathbf{W}$ and Equation (16) we deduce that $\mathbb{E}(\mathbf{R}\mathbf{R}^T)_{ij} = \delta_{ij}\mathrm{tr}\,(\mathbf{W}^T\mathbf{W})$. Similarly $\mathbb{E}(\mathbf{S}\mathbf{S}^T)_{ij} = \delta_{ij}\mathrm{tr}\,(\mathbf{V}^T\mathbf{V})$. The fact that $\mathbf{W}^T\mathbf{V} = \mathbf{0}$ makes $\mathbf{R}$ and $\mathbf{S}$ independent by Equation (16), so that:

$$\begin{aligned}
\mathbb{E}\left(\mathrm{tr}\,(\mathbf{R}\mathbf{R}^T\mathbf{P_G}\mathbf{S}\mathbf{S}^T\mathbf{P_G})\right) &= \sum_{i=1}^{n}\sum_{j=1}^{n}\sum_{k=1}^{n}\sum_{l=1}^{n}\mathbb{E}(\mathbf{R}\mathbf{R}^T)_{ij}(\mathbf{P_G})_{jk}\mathbb{E}(\mathbf{S}\mathbf{S})_{kl}(\mathbf{P_G})_{li} \\
&= \mathrm{tr}\,(\mathbf{W}^T\mathbf{W})\mathrm{tr}\,(\mathbf{V}^T\mathbf{V})\sum_{i=1}^{n}\sum_{k=1}^{n}(\mathbf{P_G})_{ik}(\mathbf{P_G})_{ki} \\
&= \mathrm{tr}\,(\mathbf{W}^T\mathbf{W})\mathrm{tr}\,(\mathbf{V}^T\mathbf{V})\mathrm{tr}\,(\mathbf{P_G}^T\mathbf{P_G})
\end{aligned}$$

as required. □

*Proof of Proposition 5.4.* During this proof, for simplicity, we write $\lambda_i = \lambda_i(\widehat{\mathbb{X}}^T\widehat{\mathbb{X}})$ and $\mathbf{L} = (\mathbb{X} - \widehat{\mathbb{X}})^T(\mathbb{X} - \widehat{\mathbb{X}})$, where Equation (34) means that as before $\mathbf{L} = \sigma^2\mathbf{X}^T\mathbf{P_G}\mathbf{X}$, with $\mathbf{P_G} = (\mathbf{I} - \frac{1}{n}\mathbf{1}_n\mathbf{1}_n^T - \mathbf{P_F})$. We condition on the $r$ PFC directions being $\mathbf{u}_i$, where $\mathbf{u}_i$ is a $\lambda_i$-eigenvector of $\widehat{\mathbb{X}}^T\widehat{\mathbb{X}}$. As in Equation (29) and [13], if the vectors $\mathbf{v}_i$ satisfy

$$\left(\mathbf{I}_p + (\widehat{\mathbb{X}}^T\widehat{\mathbb{X}} - \lambda_i\mathbf{I}_p)^+(\mathbf{L} - \mu_i\mathbf{I}_p)\right)\mathbf{v}_i = -(\widehat{\mathbb{X}}^T\widehat{\mathbb{X}} - \lambda_i\mathbf{I}_p)^+\mathbf{L}\mathbf{u}_i. \qquad (36)$$

then $\mathbf{u}_i + \mathbf{v}_i$ are $(\lambda_i + \mu_i)$-eigenvectors of $\mathbb{X}^T\mathbb{X}$, and give the $r$ PC directions. This definition of $\mathbf{v}_i$ makes $\mathbf{v}_i \perp \mathbf{u}_i$, so that $\mathbf{P}_{\mathbf{u}_i}\widehat{\mathbf{\Gamma}}_{\mathrm{PC},i} = \mathbf{u}_i$ and $(\mathbf{I}_p - \mathbf{P}_{\mathbf{u}_i})\widehat{\mathbf{\Gamma}}_{\mathrm{PC},i} = \mathbf{v}_i$. As before, we know that for all $i$:

$$\begin{aligned}
\|\mathbf{I}_p + (\widehat{\mathbb{X}}^T\widehat{\mathbb{X}} - \lambda_i\mathbf{I}_p)^+(\mathbf{L} - \mu_i\mathbf{I}_p)\| &\geq 1 - \|(\widehat{\mathbb{X}}^T\widehat{\mathbb{X}} - \lambda_i\mathbf{I}_p)^+\| \times \|(\mathbf{L} - \mu_i\mathbf{I}_p)\| \\
&\geq 1 - \frac{\|\mathbf{L}\|}{M}.
\end{aligned}$$

This means that, conditioned on $\mathcal{L}_2$,

$$\frac{\|(\mathbf{I}_p - \mathbf{P}_{\mathbf{u}_i})\widehat{\boldsymbol{\Gamma}}_{\text{PC},i}\|_F^2}{|\mathbf{u}_i|^2} \leq \frac{|(\widehat{\mathbb{X}}^T\widehat{\mathbb{X}} - \lambda_i \mathbf{I}_p)^+ \mathbf{L} \mathbf{u}_i|^2}{\|\mathbf{I}_p + (\widehat{\mathbb{X}}^T\widehat{\mathbb{X}} - \lambda_i \mathbf{I}_p)^+(\mathbf{L} - \mu_i \mathbf{I}_p)\|^2}$$

$$\leq \frac{|(\widehat{\mathbb{X}}^T\widehat{\mathbb{X}} - \lambda_i \mathbf{I}_p)^+ \mathbf{L} \mathbf{u}_i|^2}{(1 - \|\mathbf{L}\|/M)^2}.$$

Writing $\mathbf{U}$ for the matrix made up of columns $\mathbf{u}_i$, using Jensen's inequality

$$\tan^2\left(\mathbb{E}\left(\Theta(\widehat{\boldsymbol{\Gamma}}_{\text{PC}}, \mathbf{U})|\mathcal{L}_2\right)\right) \leq \mathbb{E}\left(\tan^2\Theta(\widehat{\boldsymbol{\Gamma}}_{\text{PC}}, \mathbf{U})|\mathcal{L}_2\right) \quad (37)$$

$$= \mathbb{E}\left(\left.\frac{\sum_i |(\mathbf{I}_p - \mathbf{P}_{\mathbf{U}})\widehat{\boldsymbol{\Gamma}}_{\text{PC},i}|^2}{\sum_i |\mathbf{P}_{\mathbf{U}}\widehat{\boldsymbol{\Gamma}}_{\text{PC},i}|^2}\right|\mathcal{L}_2\right)$$

$$\leq \mathbb{E}\left(\left.\frac{\sum_i |(\mathbf{I}_p - \mathbf{P}_{\mathbf{u}_i})\widehat{\boldsymbol{\Gamma}}_{\text{PC},i}|^2}{\sum_i |\mathbf{P}_{\mathbf{u}_i}\widehat{\boldsymbol{\Gamma}}_{\text{PC},i}|^2}\right|\mathcal{L}_2\right)$$

$$\leq \frac{\sum_i \mathbb{E}(((\widehat{\mathbb{X}}^T\widehat{\mathbb{X}} - \lambda_i \mathbf{I}_p)^+ \mathbf{L} \mathbf{u}_i)^2|\mathcal{L}_2)}{(1 - \|\mathbf{L}\|/M)^2 \sum_i |\mathbf{u}_i|^2}. \quad (38)$$

Using Lemma A.2 and recalling the fact that $\mathbf{L} = \sigma^2 \mathbf{X}^T \mathbf{P}_{\mathbf{G}} \mathbf{X}$, we know that the numerator of Equation (38) satisfies

$$\sigma^4 \sum_i \mathbb{E} \mathbf{u}_i^T \mathbf{X}^T \mathbf{P}_{\mathbf{G}} \mathbf{X} (\widehat{\mathbb{X}}^T\widehat{\mathbb{X}} - \lambda_i \mathbf{I}_p)^+ (\widehat{\mathbb{X}}^T\widehat{\mathbb{X}} - \lambda_i \mathbf{I}_p)^+ \mathbf{X}^T \mathbf{P}_{\mathbf{G}} \mathbf{X} \mathbf{u}_i$$

$$= \sigma^4 \sum_i \mathbb{E}\left(\text{tr}\,(\mathbf{S}_i^T \mathbf{P}_{\mathbf{G}} \mathbf{R}_i \mathbf{R}_i^T \mathbf{P}_{\mathbf{G}} \mathbf{S}_i)\right)$$

$$= \sigma^4 \sum_i \mathbb{E}\left(\text{tr}\,(\mathbf{R}_i \mathbf{R}_i^T \mathbf{P}_{\mathbf{G}} \mathbf{S}_i \mathbf{S}_i^T \mathbf{P}_{\mathbf{G}})\right)$$

$$= \sigma^4 \text{rank}\,(G) \sum_i \text{tr}\,\left((\widehat{\mathbb{X}}^T\widehat{\mathbb{X}} - \lambda_i \mathbf{I}_p)^{+2}\right) |\mathbf{u}_i|^2$$

$$\leq \frac{\sigma^4 np}{M^2} \sum_i |\mathbf{u}_i|^2,$$

where $\mathbf{R}_i = \mathbf{X}(\widehat{\mathbb{X}}^T\widehat{\mathbb{X}} - \lambda_i \mathbf{I}_p)^+$, $\mathbf{S}_i = \mathbf{X}\mathbf{u}_i$, so the conditions of Lemma A.2 apply. Substituting into Equation (38), we deduce the result. $\square$

## Acknowledgements

The author would like to thank Iain Johnstone for the reference to Davidson and Szarek's paper, and Richard Samworth for helpful advice.